\documentclass[11pt]{article}
\usepackage{latexsym, amsmath, amssymb}
\usepackage{amsmath}
\usepackage{amsfonts}
\usepackage{amssymb}
\usepackage{amscd}
\usepackage{amsthm}
\usepackage{fancyhdr}
\newtheorem{theorem}{Theorem}[section]

\newtheorem{remark}[theorem]{Remark}
\newtheorem{lemma}[theorem]{Lemma}
\newtheorem{definition}[theorem]{Definition}
\newtheorem{corollary}[theorem]{Corollary}

\def\C{\mathbb C}

\def\Q{\mathbb Q}

\newcommand{\zz}{{\mathbb Z}}

\newcommand\sH{{\mathcal H}}

\newcommand\sO{{\mathcal O}}

\newcommand\sS{{\mathcal S}}

\newcommand\sZ{{\mathcal Z}}

  \def\peen{\hbox{$ {\mathbf  P}^n$}}

  \def \tab#1{\kern #1 truein}

 \def\P#1{\hbox{${\mathbf P}^{#1}$}}

  \def\E{\hbox{${\cal E}$}}
  \def\F{\hbox{${\cal F}$}}
  \def\G{\hbox{${\cal G}$}}
  \def\A{\hbox{${{\cal A}}$}}
  \def\B{\hbox{${\cal B}$}}
  \def\C{\hbox{${\cal C}$}}

  \def\P{\hbox{${\cal P}$}}
  \def\Q{\hbox{${\cal Q}$}}

\begin{document}
\title{Monads and Rank Three Vector Bundles on Quadrics} 
\date{}
\author{Francesco Malaspina \\
 Dipartimento di Matematica Universit\`a di Torino\\
via Carlo Alberto 10, 10123 Torino, Italy\\ 
{\small\it e-mail: francesco.malaspina@unito.it}}

   \maketitle
   \def\thefootnote{}
  \baselineskip=+.5cm\begin{abstract}In this paper we give the classification of rank $3$ vector bundles without "inner" cohomology on a quadric hypersurface $\Q_n$ ($n>3$) by studying the associated monads. 
      \end{abstract}

   \section*{Introduction}

A monad on $\mathbb{P}^n$ or, more generally, on a projective variety $X$, is a complex of three vector bundles
$$0 \rightarrow \A \xrightarrow{\alpha} \B \xrightarrow{\beta} \C \rightarrow 0$$
such that $\alpha$ is injective as a map of vector bundles and $\beta$ is surjective.
Monads have been studied by Horrocks, who proved (see \cite{Ho} or \cite{BH}) that every vector bundle on $\mathbb{P}^n$  is the homology of a suitable minimal monad.
Throughout the paper we often use the Horrocks correspondence between a bundle $\E$ on $\mathbb P^n$ ($n\geq 3$) and the corresponding minimal monad
  $$0 \rightarrow \A \xrightarrow{\alpha} \B
\xrightarrow{\beta} \C \rightarrow 0,$$ where $\A$ and $\C$ are sums of line
bundles and $\B$ satisfies:
\begin{enumerate}
\item $H^1_*(\B)=H^{n-1}_*(\B)=0$ 
\item $H^i_*(\B)=H^i_*(\E)$ \ $\forall i, 1<i<n-1
$.
\end{enumerate}where $ H^i_*(\B)$ is defined as $\oplus_{k\in \mathbb(Z)} H^i(\mathbb P^n , B(k))$.\\
This correspondence holds also on a projective variety $X$ ($\dim X\geq 3$) if we fix a very ample line bundle $\sO_X(1)$. Indeed the proof of the result in (\cite{BH} proposition $3$) can be easily extended to $X$ (see \cite{Ml})). \\
Rao, Mohan Kumar and Peterson  have successfully used this tool to investigate the intermediate cohomology modules of a vector bundle on $\mathbb{P}^n$ and give cohomological splitting conditions (see \cite{KPR1}).\\
This theorem makes a strong use of monads and of Horrocks' splitting criterion which states the following:\\
Let $\E$ be a vector bundle of rank $r$ on $\mathbb P^n$, $n\geq 2$ then $\E$ splits if and only if it does not have intermediate cohomology (i.e. $H^1_*(\E)= ...=H^{n-1}_*(\E)=0$).\\
This criterion fails on  more general varieties. In fact there exist non-split vector bundles   on $X$ without intermediate cohomology. This bundles are called ACM bundles.\\
Rao, Mohan Kumar and Peterson focused on bundles without inner cohomology (i.e. $H^2_*(\E)= ...=H^{n-2}_*(\E)=0$) and showed that these bundles on $\mathbb{P}^n$ $(n>3)$ are split if the rank is small.\\

 On a quadric hypersurface $\Q_n$ the Horrocks criterion does not work, but there is a theorem that classifies all the ACM bundles (see \cite{Kn}) as  direct sums of line bundles and spinor bundles (up to a twist - for generalities about spinor bundles see \cite{Ot2}).\\

In [Ml] we improve Ottaviani's splitting criterion for vector bundles on a quadric hypersurface (see \cite{Ot1} and \cite{Ot3}) and obtain the equivalent of the result by Rao, Mohan Kumar and Peterson. Moreover we give the classification of rank $2$ bundles without inner cohomology on $\Q_n$ ($n>3$). It surprisingly  exactly agrees with the classification by Ancona, Peternell and Wisniewski of rank $2$ Fano bundles (see \cite{APW}).\\
We proved that for an  indecomposable rank $2$ bundle $\E$ on $\Q_4$ with $H_*^1(\E)\not=0$ and $H_*^2(\E)=0$, the only possible minimal monad with $\A$ and $\C$ different from zero  is (up to a twist)
$$0 \to \sO \xrightarrow{\alpha'}
\sS'(1)\oplus\sS''(1) \xrightarrow{\beta'} \sO(1) \to 0,$$
where $\sS'$ and $\sS''$ are the two spinor bundles on $\Q_4$.\\
The homology is the bundle $\sZ_4$ associated to   the disjoint union of a plane in $\Lambda$ and a plane in $\Lambda'$, the two families of planes in $\Q_4$ (see \cite{AS}).\\
The kernel $\G_4$ and the cokernel $\P_4$ (the dual) of this monad are rank $3$ bundles without inner cohomology and we have the two sequences 
\begin{equation} 0 \to \G_4 \rightarrow
\sS'(1)\oplus\sS''(1) \rightarrow \sO(1) \to 0,
\end{equation}
and
\begin{equation} 0 \to \sO \rightarrow
\sS'(1)\oplus\sS''(1) \rightarrow \P_4 \to 0.
\end{equation}
On $\Q_5$ there is only one spinor $\sS_5$ and the only possible minimal monad with $\A$ and $\C$ different from zero, for a rank $2$ bundle without inner cohomology, is (up to a twist) 
$$0 \to \sO \xrightarrow{\alpha''}
\sS(1) \xrightarrow{\beta''} \sO(1) \to 0.$$ The kernel $\G_5$ and the cokernel $\P_5$ (the dual) of the monad are rank $3$ bundles without inner cohomology and we have the two sequences 
\begin{equation} 0 \to \G_5 \rightarrow
\sS(1) \rightarrow \sO(1) \to 0,
\end{equation}
and
\begin{equation} 0 \to \sO \rightarrow
\sS(1) \rightarrow \P_5 \to 0.
\end{equation}
 The homology of the monad $\sZ_5$  is a Cayley bundle (see \cite{Ot4} for generalities on Cayley bundles).\\
The bundle $\sZ_5$ appear also in \cite{Ta} and \cite{KPR2}.\\
For $n>5$, no non-split bundle of rank $2$ on $\Q_n$ exists with\\
$H^2_*(\E)= ...=H^{n-2}_*(\E)=0$.\\
 
The main aim of the present paper is the classification of rank three bundles without inner cohomology  on $\Q_4$ by studying the associated monads.
We are able to prove that:\\
For a non-split rank $3$ bundle $\E$ on $\Q_4$ with $H_*^2(\E)=0$, the only possible minimal monads with $\A$ or $\C$ different from zero  are (up to a twist) the sequences $(1)$ and $(2)$ and
\begin{equation} 0 \to \sO \xrightarrow{\alpha}
\sS'(1)\oplus\sS''(1)\oplus \sO(a) \xrightarrow{\beta} \sO(1) \to 0,
\end{equation}
where $a$ is an integer, $\alpha=(\alpha', 0)$ and $\beta=(\beta', 0)$.\\
This means that on $\Q_4$ the only non-split rank $3$ bundles without inner cohomology are the following:\\
 the ACM bundles $\sS'\oplus\sO(a)$ and $\sS''\oplus\sO(a)$, $\G_4$, $\P_4$ and $\sZ_4\oplus\sO(a)$.\\ 
 In particular $\G_4$ and its dual are the only indecomposable rank $3$ bundles without inner cohomology on $\Q_4$.\\
 By using monads again we can also understand the behavior of rank three bundles on $Q_5$  and also on $\Q_n$, ($n >5$).\\ More precisely we can prove that:\\
 For a non-split rank $3$ bundle $\E$ on $\Q_5$ without inner cohomology, the only possible minimal monad with $\A$ or $\C$ not zero  are (up to a twist) the sequences $(3)$ and $(4)$ and
$$ 0 \to \sO \xrightarrow{\alpha}
\sS_5(1)\oplus\sO(a) \xrightarrow{\beta} \sO(1) \to 0.$$
where $a$ is an integer, $\alpha=(\alpha'', 0)$ and $\beta=(\beta'', 0)$.
This means that on $\Q_5$ the only rank $3$ bundles without inner cohomology are the following:\\
 $\G_5$, $\P_5$ and $\sZ_5\oplus\sO(a)$.\\ 
 In particular $\G_5$ and its dual are the only indecomposable rank $3$ bundles without inner cohomology on $\Q_5$.\\
   For a non-split rank $3$ bundle $\E$ on $\Q_6$ without inner cohomology, we have four possible minimal monads:
\begin{equation} 0 \to \sO \rightarrow
\sS'_6(1) \rightarrow \P'_6 \to 0,\end{equation} 
\begin{equation} 0 \to \sO \rightarrow
\sS''_6(1) \rightarrow \P''_6 \to 0,\end{equation} 
\begin{equation} 0 \to \G'_6 \rightarrow\sS'_6(1) \rightarrow \sO(1) \to 0,\end{equation} and 
\begin{equation} 0 \to \G''_6 \rightarrow
\sS''_6(1) \rightarrow \sO(1) \to 0.\end{equation} 
These sequences appear for instance in \cite{Ot2} Theorem $3.5.$\\
Therefore on $\Q_6$ the only rank $3$ bundles without inner cohomology are the following:\\
 $\G'_6$, $\G''_6$, $\P'_6$ and $\P''_6$.\\
  
   For $n>6$, no non-split bundles of rank $3$ in $\Q_n$ exist with
$H^2_*(\E)= ...=H^{n-2}_*(\E)=0$.\\
I would like to thank A. Prabhakar Rao for having introduced me into the topic and Giorgio Ottaviani for his useful comments and suggestions.\\

\section{Monads for Bundles on ACM Varieties}
In this section $X$ denotes a nonsingular subcanonical, irreducible ACM projective variety. By this we mean that $X$ has a very ample line bundle $\sO_X(1)$ such that $\omega_X \cong \sO_X(a)$ for some $a \in \mathbb{Z}$ and the embedding of $X$ by $\sO_X(1)$ is arithmetically Cohen-Macaulay. We will also assume that every line bundle on $X$ has the form $\sO_X(k), k \in \mathbb{Z}$.
 \\
If $M$ is a finitely generated graded module over the homogeneous coordinate ring of $X$, we denote by $\beta_{ij}(M)$ and $\beta_i(M)$ the graded Betti numbers and total Betti numbers of $M$. We will mainly use $\beta_{0j}(M)$ and $\beta_0(M)$ which give the number of minimal generators of $M$ in degree $j$ and the total number of minimal generators respectively.\\
We say that a bundle is non-split if it does not split as a direct sum of line bundles.\\
We say that a bundle is indecomposable if it does not split in two direct summmands.\\
\begin{definition}We will call bundle without inner cohomology a bundle $\E$ on $X$ with 
$$H^2_*(\E)= \dots =H^{n-2}_*(\E)=0,$$
where $n=dim X$.
\end{definition}
We prove a theorem about monads for rank $r$ bundles: 

\begin{theorem}\label{t2} On $X$
 of dimension $n$ with $n>3$,  any minimal monad $$ 0 \to \A \xrightarrow{\alpha} \B
\xrightarrow{\beta} \C \to 0,$$ such that $\A$ or $\C$ are not
zero, for a rank $r$ ($r\geq 2$) bundle with $H^2_*(\E)=H^{n-2}_*(\E)=H^2_*(\wedge^2 \E) =H^2_*(\wedge^2 \E^\vee)=0$, must satisfy the
following conditions:
\begin{enumerate}

\item $H^1_*(\wedge^2\B)\not=0$, $\beta_0(H^1_*(\wedge^2\B))\geq
\beta_0(H^0_*(S_2\C))$ and \\
\vskip0.01truecm 
$\beta_{0j}(H^1_*(\wedge^2\B))\geq
\beta_{0j}(H^0_*(S_2\C))$ $\forall j\in \zz$,  if $\C$ is not zero.
\vskip0.8truecm
\item $H^1_*(\wedge^2\B^{\vee})\not=0$, $ \beta_0(H^1_*(\wedge^2\B^{\vee}))\geq
\beta_0(H^0_*(S_2\A^{\vee}))$ and\\
\vskip0.01truecm $ \beta_{0j}(H^1_*(\wedge^2\B^{\vee}))\geq
\beta_{0j}(H^0_*(S_2\A^{\vee}))$ $\forall j\in \zz$, if $\A$ is not zero.
\vskip0.8truecm
 \item $H^2_*(\wedge^2\B)=H^2_*(\wedge^2\B^{\vee}) =0$.
\end{enumerate}
\begin{proof}First of all, since $X$ is ACM, the sheaf $\sO_X$ does not have intermediate
cohomology. Hence the same is true for $\A$ and $\C$ that are free $\sO_X$ sheaves.\\
Let us now  assume the existence of a minimal monad with $H^1_*(\wedge^2\B)=0$ and $\C$ not zero
$$ 0 \to \A \xrightarrow{\alpha} \B
\xrightarrow{\beta} \C \to 0.$$ Then, if we call $\G=\ker\beta$,
from the sequence
$$ 0 \to S_2\A  \to (\A \otimes
\G)\to \wedge^2 \G
   \to \wedge^2\mathcal E \to 0, $$
   we have
$$ H^2_*(\wedge^2 \G) = H^2_*(\A \otimes \G) = 0,$$ since $H^2_*(\B)=H^2_*(\G)=H^2_*(\E)=0$ and $ H^2_*(\wedge^2 \E) = 0,$.\\
Moreover, from the sequence
$$ 0 \to \wedge^2\G \to \wedge^2\B \to \B\otimes \C
   \to S_2\C \to 0, $$
 by passing to the exact
sequence of maps on cohomology groups, since $H^1_*(\wedge^2\B)=H^2_*(\wedge^2 \G)=0$ we get
$$    H^0_*(\B\otimes \C)
   \to H^0_*(S_2\C) \to 0 .$$
Now, if we call $S_X$ the coordinate ring, we can say that $H^0_*(S_2\C)$
   is a free $S_X$-module hence projective, then there exists a map $$H^0_*(\B\otimes \C)
   \leftarrow H^0_*(S_2\C)$$ and this means that $$\B\otimes \C
   \to S_2\C \to 0 $$ splits.\\
   But this map is obtained from $\beta$ as $b\otimes c\mapsto
   \beta(b)c$, so if it splits also $\beta$ has  to split and this violates the
   minimality of the monad.
We can say something stronger.\\
Recall that if $M \to N \to 0$ is a surjection of finitely generated graded $S_X$-modules, then $\beta_0(M) \geq \beta_0(N)$ and also $\beta_{0j}(M) \geq \beta_{0j}(N)$ for any $j$. Furthermore, if the inequality is strict, it means that a set of minimal generators of $M$ (in degree $j$ in the second case) can be chosen in such a way that one of generators in the set maps to zero. \\

From the sequence
$$0 \to \wedge^2\G \to \wedge^2\B \to \B\otimes \C\xrightarrow{\gamma} S_2\C \to 0, $$
    since $H^2_*(\wedge^2 \G)  = 0$, we have a surjective map
    $$    H^1_*(\wedge^2\B)
   \to H^1_*(\Gamma) \to 0 $$ where $\Gamma=\ker\gamma$, and then
   $$     \beta_0(H^1_*(\wedge^2\B))\geq \beta_0(H^1_*(\Gamma)).$$
   On the other hand we have the sequence

    $$    H^0_*(\B\otimes \C)
   \xrightarrow{\gamma} H^0_*(S_2\C)\rightarrow  H^1_*(\Gamma) \to 0; $$ so, if
    $$     \beta_0(H^1_*(\wedge^2\B))< \beta_0(H^0_*(S_2\C)),$$ also $$\beta_0(H^1_*(\Gamma))< \beta_0(H^0_*(S_2\C)),$$ and some of the generators
    of $H^0_*(S_2\C)$ must be in the image of $\gamma$.\\   
    But $\gamma$ is obtained from $\beta$ as $b\otimes c\mapsto
   \beta(b)c$, so  also  some generators of $C$ must be in the image of $\beta$ and this contradicts the
   minimality of the monad.\\
   We conclude that not just $H^1_*(\wedge^2\B)$ has to be non zero
   but also $$     \beta_0(H^1_*(\wedge^2\B))\geq \beta_0(H^0_*(S_2\C)).$$
   If we fix the degree $\j$ we have that also the map $H^0(\B\otimes \C(j))\rightarrow H^0(S_2\C(j))$ and so  we see that, $\forall \j\in\zz$
   $$     \beta_{0j}(H^1_*(\wedge^2\B))\geq \beta_{0j}(H^0_*(S_2\C)).$$
   If $\A=0$ the monad is simply $$ 0 \to \E \xrightarrow{\alpha} \B
\xrightarrow{\beta} \C \to 0,$$ and, by using the sequence $$ 0 \to \wedge^2\E \to \wedge^2\B \to \B\otimes \C\to S_2\C \to 0, $$
    since $H^2_*(\wedge^2 \E)  = 0$, we can conclude as before.\\

    Let us now  assume the existence of a monad with  $\A$ not zero $(H^1_*(\wedge^2\B^{\vee}))=0$, we use the dual sequences.\\
    From $$ 0 \to S_2\C^{\vee}  \to (\C^{\vee} \otimes
\B^{\vee})\to \wedge^2 \B^{\vee}
   \to \wedge^2\G^{\vee} \to 0, $$
   we have
$ H^1_*(\wedge^2 \G^{\vee})\cong H^1_*(\wedge^2 \B^{\vee})$.\\
Moreover, from the sequence
$$ 0 \to \wedge^2\E^{\vee} \to \wedge^2\G^{\vee} \to \G^{\vee}\otimes \A^{\vee}
   \to S_2\A^{\vee} \to 0, $$
 by passing to the exact
sequence of maps on cohomology groups, since $H^2_*(\wedge^2\E^{\vee})=H^1_*(\wedge^2 \G^{\vee})=0$ we get
$$    H^0_*(\G^{\vee}\otimes \A^{\vee})
   \to H^0_*(S_2\A^{\vee}) \to 0 ,$$ and this violates the
   minimality of the monad as before.\\
   We can also conclude that  $$     \beta_0(H^1_*(\wedge^2\B^{\vee}))=\beta_0(H^1_*(\wedge^2\G^{\vee}))\geq \beta_0(H^0_*(S_2\A^{\vee})).$$
   If we fix the degree $\j$ we have that also the map $H^0(\G^{\vee}\otimes \A^{\vee}(j))\rightarrow H^0(S_2\A^{\vee}(j))$ and so we see that, $\forall \j\in\zz$
   $$     \beta_{0j}(H^1_*(\wedge^2\B^\vee))\geq \beta_{0j}(H^0_*(S_2\A^\vee)).$$
   If $\C=0$ the monad is simply $$ 0 \to \A \xrightarrow{\alpha} \B
\xrightarrow{\beta} \E \to 0,$$ and, by using the sequence $$ 0 \to \wedge^2\E^{\vee} \to \wedge^2\B^{\vee} \to \B^{\vee}\otimes \A^{\vee}
   \to S_2\A^{\vee} \to 0, $$
    since $H^2_*(\wedge^2 \E^\vee)  = 0$, we can conclude as before.\\

   The third condition
comes from the sequences $$ 0 \to \wedge^2\G \to \wedge^2\B \to
\B\otimes \C
   \to S_2\C \to 0, $$ and $$ 0 \to S_2\C^{\vee}  \to (\C^{\vee} \otimes
\B^{\vee})\to \wedge^2 \B^{\vee}
   \to \wedge^2\G^{\vee} \to 0, $$ since $H^2_*(\wedge^2\G)=H^2_*(
\B\otimes \C)=H^2_*((\C^{\vee} \otimes
\B^{\vee})=H^2_*(\wedge^2\G^{\vee})=0$.
   \end{proof}
   \end{theorem}
   \begin{remark}If $r=2$ we have ([ML]$1.6$).
   \end{remark}
   \begin{remark} If $r=3$ we don't need the hypothesis $H^2_*(\wedge^2 \E) =H^2_*(\wedge^2 \E^\vee)=0$ because $H^2_*(\wedge^2 \E) = H^{n-2}_*(\E) = 0.$
   \end{remark}
   \begin{remark}On $\mathcal{P}^n$ we can say the following:\\
   Let $\E$ be a bundle without inner cohomology such that $H^2_*(\wedge^2 \E) =H^2_*(\wedge^2 \E^\vee)=0$, then $\E$ splits.\\
   In fact if $\E$ does not split the associated minimal monad has $\A$ or $\C$ different to zero. Since $H^2_*(\E)= ...=H^{n-2}_*(\E)=0$, the bundles $\B$ does not have intermediate cohomology and hence it splits. In particular $H^1_*(\wedge^2 \B) =0$.  By the above theorem this is a contradiction.\\
   So the hypothesis $H^2_*(\wedge^2 \E) =H^2_*(\wedge^2 \E^\vee)=0$ avoid the limitation of the rank in the Kumar, Peterson and Rao theorem  (see \cite{KPR1}).
   \end{remark}
   We need also the following lemma:\\
   \begin{lemma} Let $\E$ be a rank $2$ on $X$. If 
$$ 0 \to \A \xrightarrow{\alpha'}
\B \xrightarrow{\beta'} \C \to 0,$$  is a minimal monad for $\E$, then  
$$ 0 \to \A \xrightarrow{\alpha}
\B\oplus \sO(a) \xrightarrow{\beta} \C \to 0,$$  where $\alpha=(\alpha', 0)$ and $\beta=(\beta', 0)$, is a minimal monad for $\E\oplus \sO(a)$.
\end{lemma}

   \section{Rank $3$ Bundles without Inner\\ Cohomology}

   We want now apply these results in order to classify the rank $3$ bundles without inner cohomology on $\Q_4$:
   \begin{theorem}For a  non-split rank $3$ bundle $\E$ on $\Q_4$ with\\ $H_*^2(\E)=0$, the only possible minimal monads with $\A$ or $\C$ different from zero  are (up to a twist) the sequences $(1)$ and $(2)$ and
\begin{equation} 0 \to \sO \rightarrow
\sS'(1)\oplus\sS''(1)\oplus \sO(a) \rightarrow \sO(1) \to 0,
\end{equation}
where $a$ is an integer, $\alpha=(\alpha', 0)$ and $\beta=(\beta', 0)$.
\begin{proof}First of all consider a minimal  monad for $\E$,
$$ 0 \to \A \xrightarrow{\alpha} \B
\xrightarrow{\beta} \C \to 0,$$ Since by construction, $\B$  is an ACM bundle on $\Q_4$,
then it has to be isomorphic to a direct sum of line bundles and
spinor bundles twisted by some $\sO(t)$.\\
  Since $\B$ cannot be split without violating part $1$ of (Theorem \ref{t2}) which states that $$H^1_*(\wedge^2\B)\not=0.$$ Hence at least a spinor bundle must appear in its decomposition.\\
If just one copy of $\sS'$ or one copy of $\sS''$ appears in $\B$, since $$\textrm{rank $\sS''$ $=$ rank $\sS'$ $=2$},$$ and then $\wedge^2 \sS'$ and $\wedge^2\sS''$ are line bundles, also the
bundle $\wedge^2\B$ is ACM and again the condition
$$H^1_*(\wedge^2\B)\not=0,$$ in (Theorem \ref{t2}), is not satisfied.\\
 If it appears more than one copy of $\sS'$ or more than one copy of $\sS''$ appears in $\B$, in the
bundle $\wedge^2\B$ it appears $(\sS'\otimes\sS')(t)$ or
$(\sS''\otimes\sS'')(t)$ appears and, since
$$H^2_*(\sS'\otimes\sS')=H^2_*(\sS''\otimes\sS'')=\mathbb
C,$$  the condition
$$H^2_*(\wedge^2\B)=0$$ in (Theorem (\ref{t2})), fails to be satisfied. So $\B$ must contain both $\sS'$ and $\sS''$ 
with some twist and only one copy of each. We can conclude that $\B$  has to be of the form $$ 
 (\bigoplus_i\sO(a_i))\oplus ( \sS'(b))\oplus (\sS''(c)).$$ Let us notice
furthermore that if $H_*^1(\E)$ has more than $1$ generator, rank $
\C>1$ and $H_*^0( S_2\C)$ has at least $3$ generators.\\
But
$$H_*^1(\wedge^2\B)\backsimeq H_*^1(\sS'\otimes\sS'')=\mathbb C$$ has just $1$ generator and
this is a
contradiction because by (Theorem (\ref{t2})) $$ \beta_0(H^1_*(\wedge^2\B))\geq
\beta_0(H^0_*(S_2\C)).$$ 
This means that rank $\C$ $=1$ or $=0$ .\\
Similarly, looking at dual sequence, we have that also rank $\A$ must be $1$.\\
If $\C=0$, we have the  minimal monad 
$$0 \to \sO \rightarrow
\sS'(l)\oplus\sS''(m) \rightarrow \E \to 0.$$ By computing $c_4(\sS'(l)\oplus\sS''(m))$ as in (\cite{Ml} Theorem $3.1$) we see that $l$ and $m$ must be both equal to $1$ and we have the monad ($2$).\\
If $\A=0$ we see in the same way that we have the monad ($1$).\\
At this point the only possible monads with $\A$ and $\C$ not zero, are like
$$ 0 \to \sO \xrightarrow{\alpha}
\sO(a)\oplus\sS'(1+b)\oplus\sS''(1+c) \xrightarrow{\beta} \sO(d) \to 0.$$ where $a, b$ and $c$
are integer numbers.\\
Now, since $$\beta_{0j}(H^1_*(\wedge^2\B))\geq
\beta_{0j}(H^0_*(S_2\C))$$ $\forall j\in \zz$ we see that $2+b+c=2d$.\\
 Let assume that $b\leq 0$ then by the sequence $$ 0 \to \sS'' \rightarrow
\sO^4 \rightarrow \sS'(1) \to 0,$$ (see \cite{Ot2}) we see that $\sS'(1+b)$ does not have global section.\\
Moreover $\sO(a)\oplus\sS''(1+c)$ does not have nowhere vanishing section. In fact a section of $\sO(a)$ has zero locus of dimension $4$ (if it is the zero map) or $3$. If $a=0$ it could be a scalar different to zero but this is against our assumption of minimality. Since the zero locus of a section of $\sS''(1+c)$ has dimension at least $2$, we conclude that the zero locus of a section of $\sO(a)\oplus\sS''(1+c)$ must be not empty.\\
This means that the map $\alpha$ cannot be injective.\\
 If $c\leq 0$ we have the same contradiction.\\
 We have, hence, that $b$ and $c$ must be positives.
Let us consider now the dual monad twisted by $d$ $$ 0 \to \sO \xrightarrow{\beta^\vee}
\sO(d-a)\oplus\sS'(d-b)\oplus\sS''(d-c) \xrightarrow{\alpha^\vee} \sO(d) \to 0.$$
By the argument above we have that $d-b\geq 1$ and $d-c\geq 1$; but, since $2=d-b+d-c$, it follows that $b=c$ and $d=b+1$.\\
We have the map $$\beta: \sO(a) \oplus \sS'(1+b)\oplus \sS''(1+b) \rightarrow \sO(1+b).$$
Let us consider the restriction $$\beta' : \sS'(1+b)\oplus \sS''(1+b) \rightarrow \sO(1+b).$$ We know (by \cite{Ml}) that in general, we can find a surjective map $$\gamma: \sS'(1+b)\oplus \sS''(1+b) \rightarrow \sO(1+b)$$ since we have some standard rank two bundles obtained from a monad on $\Q_4$. Hence the map

$\gamma^\vee$ gives a nowhere vanishing section of $\sS'^\vee \oplus \sS''^\vee$, which thus has fourth Chern class $0$. Hence in particular,  some other map like

$\beta'^\vee$ must give a section which is either nowhere vanishing, or which vanishes on a locus of dimension $\geq 1$. (It cannot vanish on a non empty zero dimensional set).  However, if $\beta'^\vee$ vanishes on a locus of dimension $\geq 1$, then $\beta^\vee$ itself cannot give a nowhere vanishing section since the map $\sO  \rightarrow \sO(-a+1+b)$ is either zero or defines a hypersurface, by minimality.

Therefore $\beta'$ must be surjective (like a standard map $\gamma$).

By an easy computation we have the following claim:\\ If $\E$ is a rank two bundle on $\Q_4$ with monad  $$0 \rightarrow \sO \rightarrow \sS'(1) \oplus \sS''(1) \rightarrow \sO(1) \to 0,$$ then $H^1(E(-1)) = k$ and $H^1(E(t)) = 0$ for every $t\not= -1$.\\

Hence on the level of global sections, $\beta'$ is surjective onto $\sO(1+b)$, except that the element $1$ in degree $(-1-b)$ is not in the image. By minimality, $1$ cannot be in the image of $\sO(a)$. Hence $\sO(a)$ maps by $\beta$ to the image of $\beta'$ i.e. there exists $l\in \sS'(1+b)\oplus\sS''(1+b)$ such that $\beta(1,0)=\beta'(l)$ . Therefore, after a change of basis, we may assume that $\sO(a)$ maps to zero. In fact, if we consider a map $$\delta:  \sO(a)\oplus\sS'(1+b)\oplus\sS''(1+b) \rightarrow \sO(a)\oplus\sS'(1+b)\oplus\sS''(1+b)$$ sending $(1,0)$ in $(1,-l)$, we have that $$\beta(\delta(1,0))=\beta(1,-l)=\beta'(l)-\beta'(l)=0.$$
We have at this point the monad $$ 0 \to \sO \xrightarrow{(h, \alpha')}
\sO(a)\oplus\sS'(1+b)\oplus\sS''(1+b) \xrightarrow{(0, \beta')} \sO(1+b) \to 0.$$
We want to prove that $h$ must be the zero map and $b$ must be zero.\\
If $a\leq 0$, clearly $h=0$ and $\alpha'$ is injective if and only if $b=0$.\\
If $a>0$ we consider the kernel of $\beta$ $\sO(a)\oplus\G_4(b)$ and, from the exact sequence $$ 0 \to \sO \rightarrow
\sO(a)\oplus\G_4(b) \rightarrow \E \to 0,$$ we see that $c_4(\sO(a)\oplus\G_4(b))$ must be zero. But from $$0 \to \G_4(b) \rightarrow
\sS'(1+b)\oplus\sS''(1+b) \rightarrow \sO(1+b) \to 0,$$	we see that $c_3(\G_4(b))=c_4(\sS'(1+b)\oplus\sS''(1+b))*c_1(\sO(1+b))^{-1}=2(1+b+b^2)(b+b^2)(1+b)^{-1}$ and so $c_4(\sO(a)\oplus\G_4(b))=c_1(\sO(a))*c_3(\G_4(b))=0$ if and only if $b=0$.

\end{proof}
\end{theorem}
\begin{remark} On $\Q_4$ the only rank $3$ bundles without inner cohomology are the ACM bundles, $\G_4$, $\P_4$ and $\sZ_4\oplus\sO(a)$.
\end{remark}
\begin{corollary} In higher dimension we have:
\begin{enumerate}
\item For a non-split rank $3$ bundle $\E$ on $\Q_5$ without inner cohomology, the only possible minimal monad with $\A$ or $\C$ not zero  are (up to a twist) the sequences $(3)$ and $(4)$ and
$$ 0 \to \sO \rightarrow
\sS_5(1)\oplus\sO(a) \rightarrow \sO(1) \to 0.$$
where $a$ is an integer, $\alpha=(\alpha'', 0)$ and $\beta=(\beta'', 0)$.
\item For a non-split rank $3$ bundle $\E$ on $\Q_6$ without inner cohomology, the only possible minimal monad with $\A$ or $\C$ not zero  are (up to a twist) the sequences $(6)$, $(7)$, $(8)$ and $(9)$.
 \item For $n>6$, no non-split bundle of rank $3$ in $\Q_n$ exist with\\
$H^2_*(\E)= ...=H^{n-2}_*(\E)=0$.
\end{enumerate}

\begin{proof}
First of all let us notice that for $n>4$ there is not non-split ACM rank $3$ bundles
since the spinor bundles have rank greater than $3$.\\
Let us assume
then that $H^1_*(\E)\not=0$ or $H^{n-1}_*(\E)\not=0$ and let us see how many monads it is possible
to find:
\begin{enumerate}
\item
In a minimal monad for $\E$ on $\Q_5$,
$$ 0 \to \A \xrightarrow{\alpha} \B
\xrightarrow{\beta} \C \to 0,$$ $\B$  is an ACM bundle on $\Q_5$;
then it has to be isomorphic to a direct sum of line bundles and
spinor bundles twisted by some $\sO (t)$,

Moreover, since $H_*^2(\E)=0$ and $H_*^3(\E)=0$,
$\E_{|\Q_4}=\F$ is a bundle with $H^2_*(\F)=0$ and by (\cite{Ml} Lemma $1.2$) his minimal
monad is just the restriction of the minimal  monad for $\E$ $$ 0 \to \A \xrightarrow{\alpha} \B
\xrightarrow{\beta} \C \to 0.$$
For the theorem above, hence, this monad must be
$$ 0 \to \sO \rightarrow
\sS'(1)\oplus\sS''(1)\oplus\sO(a) \rightarrow \sO(1) \to 0.$$

Now, since $$\sS_{5_{|\Q_4}}\backsimeq \sS'\oplus\sS'',$$ the only bundle of the form 
$$(\bigoplus_i\sO(a_i))\oplus ( \bigoplus_j\sS_5(b_j))$$ having $\sS'(1)\oplus\sS''(1)\oplus\sO(a)$ as restriction on $\Q_4$ is $\sS_5(1)\oplus\sO(a)$ and then if $\A$ and $\C$ are different to zero the claimed
monad $$ 0 \to \sO \xrightarrow{\alpha} \sS_5(1)\oplus\sO(a)
\xrightarrow{\beta} \sO(1) \to 0$$ where $\alpha=(\alpha'', 0)$ and $\beta=(\beta'', 0)$, is the only possible.\\
If $\A=0$, we have the monad ($3$) and if $\C=0$, we have the monad ($4$).\\

 \item
In $\Q_6$ we use the same argument. Let us consider a minimal monad for $\E$.\\
If $\A$ and $\C$ are not zero the restriction of the monad on $\Q_5$ must be $$ 0 \to \sO \xrightarrow{\alpha} \sS_5(1)\oplus\sO(a)
\xrightarrow{\beta} \sO(1) \to 0.$$ Since
$\sS'_{6_{|\Q_5}}\backsimeq \sS_5$
and also
$\sS''_{6_{|\Q_5}}\backsimeq \sS_5,$
we have two possible minimal monads:
$$ 0 \to \sO \rightarrow
\sS'_6(1)\oplus\sO(a) \rightarrow \sO(1) \to 0$$ and $$ 0 \to \sO \rightarrow
\sS''_6(1)\oplus\sO(a) \rightarrow \sO(1) \to 0,$$ where the maps are of the form $(\gamma, 0)$.\\
In both the sequences  the homology is a bundle $\F\oplus\sO(a)$ where $\F$ is a rank $2$ bundle without inner cohomology that by (\cite{Ml} Cor. $3.4$) cannot exist, so they cannot be the monads of a rank $3$ bundles.\\
If $\A$ or $\C$ are  zero the restriction of the minimal monad on $\Q_5$ must be the minimal monad ($3$) or the minimal monad $(4)$. We have four possible minimal monads:
$$ 0 \to \sO \rightarrow
\sS'_6(1) \rightarrow \E \to 0,$$ $$ 0 \to \sO \rightarrow
\sS''_6(1) \rightarrow \E \to 0,$$ 
$$ 0 \to \E^{\vee} \rightarrow\sS'_6(1) \rightarrow \sO(1) \to 0,$$ and $$ 0 \to \E^{\vee} \rightarrow
\sS''_6(1) \rightarrow \sO(1) \to 0.$$ 
These are precisely the sequences $(6)$, $(7)$, $(8)$ and $(9)$.
\item
 Let us consider a minimal monad for bundle without inner cohomology $\E$ on $\Q_7$:
 $$ 0 \to \A \rightarrow
\B \rightarrow \C \to 0.$$ 
$\B$ must be not split and ACM. Since $\sS_{7_{|\Q_5}}\backsimeq \sS'_6\oplus\sS''_6,$
the restriction of the monad on $\Q_6$ cannot be one of the sequence $(6)$, $(7)$, $(8)$ and $(9)$.\\
We can conclude that  no non-split bundle of rank $3$ in
$\Q_7$ exists without inner cohomology.\\
Clearly also in
higher dimension it is not possible to find any rank $3$ bundle without inner cohomology.
\end{enumerate}
 
\end{proof}
\end{corollary}
\begin{remark} On $\Q_n$ ($n>3$) the only rank $3$ bundles without inner cohomology are the following:\begin{enumerate}
\item for $n=4$, the ACM bundles $\sS'\oplus\sO(a)$ and $\sS''\oplus\sO(a)$, $\G_4$, $\P_4$ and $\sZ_4\oplus\sO(a)$.\\
\item For $n=5$, $\G_5$, $\P_5$ and $\sZ_5\oplus\sO(a)$.\\
\item For $n=6$, $\G'_6$, $\G''_6$, $\P'_6$ and $\P''_6$.\end{enumerate}
\end{remark}
\begin{remark} If we consider rank $r$ ($r\geq 4$) without inner cohomology  we cannot have such a simple classification on $\Q_n$ ($n\geq 4$).\\
In fact  let $\sH$ be any ACM bundle of rank $r$ ($r > 4$) on $\Q_4$. The generic  map $$0 \to \sO^{r-4} \xrightarrow{\alpha}
\sH$$ is injective, so the cokernel of $\alpha$ is a rank $4$ bundle without inner cohomology.\\
This means there are many bundles without inner cohomology of rank $r$ ($r\geq 4$) on $\Q_n$ ($n\geq 4$).
\end{remark}
\bibliographystyle{amsplain}

\end{document}